\documentclass[12pt]{article}

\usepackage{graphicx,amsfonts,latexsym,a4wide}

\usepackage{algorithmic}
\usepackage{algorithm}
\usepackage{upquote}

\usepackage{amssymb}
\usepackage{amsmath}
\usepackage{amsthm}
\usepackage{amsfonts}

\newcommand{\be}{\begin{equation}}
\newcommand{\ee}{\end{equation}}
\newcommand{\bean}{\begin{eqnarray*}}
\newcommand{\eean}{\end{eqnarray*}}
\newcommand{\bea}{\begin{eqnarray}}
\newcommand{\eea}{\end{eqnarray}}
\newcommand{\ba}{\begin{array}}
\newcommand{\ea}{\end{array}}

\newcommand{\ul}{\underline}

\newcommand{\cB}{{\cal B}}
\newcommand{\cC}{{\cal C}}

\newcommand{\vol}{\mbox{\bf Vol}}

\newcommand{\Prob}{{\sf P}}

\newcommand{\SR}{\mathbb{R}}

\newcommand{\ME}{{\sf E\,}}

\newtheorem{theorem}{Theorem} 
\newtheorem{proposition}{Proposition}
\newtheorem{lemma}{Lemma}
\newtheorem{corollary}{Corollary}

\newtheorem{definition}{Definition}
\newtheorem{fig}{Fig.}
\newtheorem{example}{Example}

\newcommand{\bth}{\begin{theorem}}
\newcommand{\ethm}{\end{theorem}}
\newcommand{\bprop}{\begin{proposition}}
\newcommand{\eprop}{\end{proposition}}
\newcommand{\blem}{\begin{lemma}}
\newcommand{\elem}{\end{lemma}}
\newcommand{\bdf}{\begin{definition}}
\newcommand{\edf}{\end{definition}}
\newcommand{\bfig}{\begin{fig}}
\newcommand{\efig}{\end{fig}}
\newcommand{\bex}{\begin{example}}
\newcommand{\eex}{\end{example}}

\begin{document}

\begin{center}
\setlength{\baselineskip}{22pt}
{\Large \bf Why does Monte Carlo Fail to Work Properly \\in High-Dimensional Optimization Problems?}

\vskip .2in
\setlength{\baselineskip}{20pt}
{\large \bf Boris Polyak,  Pavel Shcherbakov}

Institute for Control Science, RAS, Moscow, Russia\\
E-mails: boris@ipu.ru;  cavour118@mail.ru


July 5, 2016


\end{center}

\begin{abstract}
The paper proposes an answer to the question formulated in the title.
\end{abstract}

\setlength{\baselineskip}{17pt}
\section{Introduction}

After the invention of the Monte Carlo (MC) paradigm by S.~Ulam in the late 1940s, it has become extremely popular
in numerous application areas such as physics, biology, economics, social sciences, and other areas.
As far as mathematics is concerned,
Monte Carlo methods showed themselves exceptionally efficient in the simulation of various probability distributions,
numerical integration, estimation of the mean values of the parameters, etc. \cite{MChandbook,TeCaDa}.
The salient feature of this approach to solution of various problems of this sort is that ``often,'' it is dimension-free
in the sense that, given~$N$ samples, the accuracy of the result does not depend on the dimension of the problem.

On the other hand, applications of the MC paradigm in the area of optimization are not that triumphant.
In this regard, problems of global optimization deserve special attention.
As explained in~\cite{ZhZh} (see beginning of Chapter 1.2),
``\emph{In global optimization, randomness can appear in several ways.
The main three are: (i)~the evaluations of the objective function are corrupted by random errors;
(ii)~the points $x_i$ are chosen on the base of random rules, and
(iii)~the assumptions about the objective function are probabilistic.}''
Pertinent to the exposition of this paper is only case (ii).
Monte Carlo is the simplest, brute force example of randomness-based methods (in \cite{ZhZh} it is referred to as ``Pure Random Search'').
With this method, one samples points uniformly in the feasible domain, computes the values of the objective function,
and picks the record value as the output.

Of course, there are dozens of more sophisticated stochastic methods such as
simulated annealing, genetic algorithms, evolutionary algorithms, etc.; e.g., see \cite{ZhZh,Pardalos,Simon,Goldberg,Wets}
for an incomplete list of relevant references.
However most of these methods are heuristic in nature; often, they lack rigorous justification, and the computational efficiency is questionable.
Moreover, there exist pessimistic results on ``insolvability of global optimization problems.''
This phenomenon has first been observed as early as in the monograph~\cite{NemYud} by A.~Nemirovskii and D.~Yudin,
both in the deterministic and stochastic optimization setups (see Theorem, Section 1.6 in \cite{NemYud}).
Specifically, the authors of~\cite{NemYud} considered the minimax approach to the minimization of the class of Lipschitz functions
and proved that, no matter what is the optimization method, it is possible to construct a problem which will require
exponential (in the dimension) number of function evaluations. The ``same'' number of samples is required for the simplest
MC method. Similar results can be found in~\cite{Nesterov}, Theorem 1.1.2, where the construction of ``bad'' problems
is exhibited.
Below we present another example of such problems (with very simple objective functions, close to linear ones)
which are very hard to optimize.
Concluding this brief survey, we see that any advanced method of global optimization cannot outperform Monte Carlo
when optimizing ``bad'' functions.

This explains our interest in the MC approach as applied to the optimization setup.
In spite of the pessimistic results above, there might be a belief that, if Monte Carlo is applied to
a ``good'' optimization problem (e.g., a convex one), the results would not be so disastrous.
The goal of the present paper is to clarify the situation.
We examine the ``best'' optimization problems (the minimization of a linear function on a ball or on a cube) and
estimate the accuracy of the Monte Carlo method. Unfortunately, the dependence on the dimension remains exponential, and practical
solution of these simplest problems via such an approach is impossible for high dimensions.

The second goal of the paper is to exhibit the same situation with multiobjective optimization~\cite{Deb}.
We treat methods for the Pareto set discovery via the Monte Carlo technique and estimate their accuracy, which happens to be poor for large dimensions.
These results are instructive for multiobjective optimization, because there exist many methods based on a similar approach
(with regular grids instead of random grids), see \cite{PSI,StatnikovEtAl-JOTA}.

An intuitive explanation of the effects under consideration can rely on the geometric nature of multidimensional spaces.
Numerous facts of this sort can be found in Chapter~2 of the book~\cite{HK}, which is available in the internet.
The titles of sections in~\cite{HK} are very impressive: ``The Volume is Near the Equator,'' ``The Volume is in a Narrow Annulus,''
``The Surface Area is Near the Equator.''
Some of the results in the present paper clarify these statements by providing rigorous closed-form estimates for the
minimum number of random points in the ball-shaped sets required to assess, with a given probability, the optimum of a linear function with given accuracy.
These estimates are based on our previous results on the properties of the uniform distribution over a ball,~\cite{spher} (see Section 2.2).

As far as the geometry of many-dimensional spaces is concerned, the highly advanced monograph~\cite{Milman} is worth mentioning;
it deals with the geometrical structure of finite dimensional normed spaces, as the dimension grows to infinity, and presents numerous
deep mathematical results in the area.

The rest of the paper is organized as follows. In Section~\ref{S:statement}, we propose a motivating example,
formulate the two optimization problems,
scalar and multiobjective, considered within the Monte Carlo setup, and present two known theorems on the uniform distribution over the $l_2$-norm ball.
These theorems will be used in Section~\ref{S:ball} to derive new results
related to the two optimization problems of interest for the case of the $l_2$-norm ball.
Section~\ref{S:box} deals with the scalar optimization problem for the case where~$Q$ is a box; use of various deterministic grids are also discussed.
Brief concluding remarks are given in the last section.

A preliminary version of this paper is~\cite{MC_Gran}. Several refinements are performed in the current text;
first, we changed the overall structure of the exposition; then, we provide a new result on the probability of the empirical maximum
of a linear function on a ball (Section~\ref{ssec:scalar}),
next, we add a result on the expected value (end of Section~\ref{ssec:multiobj}), present closed-form results for the $l_\infty$- and
$l_1$-norm balls (Section~\ref{ssec:boxMC}), discuss deterministic grids over a box (Section~\ref{ssec:lptau}), and accompany these new results with
numerical illustrations.
Finally, the introduction section and the bibliography list are considerably modified and extended and various typos and inaccuracies are corrected.

\section{Statement of the Problem}
\label{S:statement}

In this section, we propose a motivation for the research performed in this paper, formulate the problems of interest,
and present two known facts which form the basis for deriving the new results in Section~\ref{S:ball}.

\subsection{A Motivating Example}
To motivate our interest in the subject of this paper, we present a simple example showing failure of stochastic global optimization
methods in high-dimensional spaces.
This example is constructed along the lines suggested in~\cite{NemYud} (also, see~\cite{Nesterov}, Theorem 1.1.2)
and is closely related to one of the central problems discussed below, the minimization of a linear function over a ball in ${\mathbb R}^n$.

Consider an unknown vector $c\in {\mathbb R}^n$, $\|c||=1$, and the function
$$
f(x)=\min \Bigl\{99-c^\top x, \,\bigl(c^\top x-99\bigr)/398\Bigr\}
$$
to be minimized over the Euclidean ball $Q\subset {\mathbb R}^n$ of radius $r=100$ centered at the origin.
Obviously, the function has one local minimum $x_1=-100c$, with the function value $f_1=-0.5$,
and one global minimum $x^*=100c$, with the function value $f^*=-1$.
The objective function is Lipschitz with Lipschitz constant equal to $1$, and $\max f(x) - \min f(x)=1$.

Any standard (not problem-oriented) version of stochastic global search (such as multistart, simulated annealing, etc.)
 will miss the domain of attraction of the global minimum with probability
$1-V^1/V^0$, where $V^0$ is the volume of the ball $Q$, and $V^1$ is the volume of the set $C = \{x\in Q\colon c^\top x\ge 99\}$.
In other words, the probability of success is equal to
$$
{\sf P}= \frac{V^1}{V^0} = \frac{1}{2}I\Bigl(\dfrac{2rh-h^2}{r^2};\dfrac{n+1}{2},\dfrac{1}{2}\Bigr),
$$
where $I(x; a, b)$ is the regularized incomplete beta function with parameters~$a$ and~$b$ (for use of this function, also see
Theorem~\ref{th:ball_maxlin} in Section~\ref{S:ball}), and $h$ is the height of the spherical cap~$C$; in this example, $h=1$.
This probability quickly goes to zero as the dimension of the problem grows; say, for $n=15$ it is of the order of $10^{-15}$.
Hence, any ``advanced'' method of global optimization will find the minimum with relative error not less $50\%$;
moreover, such methods are clearly seen to be no better than a straightforward Monte Carlo sampling. The same is true if our
goal is to estimate the minimal value of the function $f^*$ (not the minimum point $x^*$). Various methods based on ordered
statistics of sample values (see Section 2.3 in \cite{ZhZh}) fail to reach the set $C$ with high probability, so that the prediction will
be close to $f_1=-0.5$ instead of $f^*=-1$.

These observations motivate our interest in the analysis of the performance of the MC schemes in optimization.

\subsection{The Two Optimization Settings}
\label{two_settings}

Let~$Q\subset\mathbb{R}^n$ denote a unit ball in one or another norm and let $\left.\xi^{(i)}\right|_1^N = \bigl\{\xi^{(1)},\dots,\xi^{(N)}\bigr\}$
be a multisample of size~$N$ from the uniform distribution $\xi\sim{\cal U}(Q)$.

We are targeted at solving the problems of the following sort.

\vskip .1in
{\bf I. Scalar Optimization:} Given the scalar-valued linear function
\begin{equation}
\label{linfun}
g(x) = c^\top x, \quad c\in\mathbb{R}^n,
\end{equation}
defined on the unit ball~$Q\subset\mathbb{R}^n$, estimate its maximum value from the multisample.

More specifically, let $\eta^*$ be the true maximum of $g(x)$ on~$Q$ and let
\begin{equation}
\label{empir_max}
\eta = \max\{g^{(1)}, \dots, g^{(N)}\}, \qquad g^{(i)}= g(\xi^{(i)}),\quad i=1,\dots,N,
\end{equation}
be the empirical maximum; we say that $\eta$ approximates~$\eta^*$ \emph{with accuracy at least~$\delta$} if
$$
\frac{\eta^* - \eta}{\eta^*} \,\leq\, \delta.
$$

Then the problem is: \emph{Given a probability level~$p\in (0,\, 1)$ and accuracy~$\delta\in (0,\,1)$,
determine the minimal length~$N_{\min}$ of the multisample such that, with probability at least~$p$,
the accuracy of approximation is at least~$\delta$ (i.e.,  with high probability, the empirical maximum nicely evaluates the true one).}

These problems  are the subject of discussion in Sections \ref{ssec:scalar} and~\ref{ssec:boxMC}.

\vskip ,1in
{\bf II. Multiobjective Optimization:} Consider now $1<m<n$ scalar functions $g_j(x),\;j=1,\dots,m$, and the image of~$Q$ under these mappings.
The problem is to ``characterize'' the boundary of the image set $g(Q)\subset\mathbb{R}^m$ via the multisample $\left.\xi^{(i)}\right|_1^N$ from~$Q$.

In rough terms, the problem is: \emph{Determine the minimal sample size~$N_{\min}$ which guarantees, with high probability, that the image of at least
one sample fall close to the boundary of~$g(Q)$}.

For the case where~$Q$ is the Euclidean ball, the mappings $g_j(x)$ are linear, and $m=2$, this problem is discussed in Section~\ref{ssec:multiobj};
various statistics (such as the cumulative distribution function, mathematical expectation, mode) of a specific random variable associated with
image points are evaluated.

\subsection{Supporting Material}

The results presented in Section~\ref{S:ball} below are based on the following two facts established in~\cite{spher};
they relate to the probability distribution of a specific linear or quadratic function of the random vector uniformly distributed on the Euclidean ball.

{\bf Fact 1} \cite{spher}.
\emph{Let the random vector $\xi\in\mathbb{R}^n$ be uniformly distributed on the unit Euclidean ball~$Q\subset\mathbb{R}^n$.
Assume that a matrix~$A\in\mathbb{R}^{m\times n}$ has rank $m\leq n$. Then the random variable
$$
\rho\doteq \Bigl( (AA^\top)^{-1}A\xi,\,A\xi \Bigr)
$$
has the beta distribution~${\cal B}(\frac{m}{2},\,\frac{n-m}{2}+1)$ with probability density function
\begin{equation}
\label{bet_distr}
f_\rho(x) =\left\{
\begin{array}{cl}
\displaystyle{\frac{\Gamma(\frac{n}{2}+1)}{\Gamma(\frac{m}{2})\Gamma(\frac{n-m}{2}+1)}\,
x^{\frac{m}{2}-1}(1-x)^{\frac{n-m}{2}}}  & \mbox{~~~for~} 0 \leq x \leq 1,  \\
\displaystyle{~~~~0} & \mbox{~~~otherwise},
\end{array}
           \right.
\end{equation}
where~$\Gamma(\cdot)$ is the Euler gamma function.}%

\vskip .1in
Alternatively, the numerical coefficient in~\eqref{bet_distr} writes
$\frac{\Gamma(\frac{n}{2}+1)}{\Gamma(\frac{m}{2})\Gamma(\frac{n-m}{2}+1)} = 1/B(\frac{m}{2},\frac{n-m}{2}+1)$,
where $B(\cdot,\cdot)$ is the beta function.

\vskip .1in
The second fact is an asymptotic counterpart of Fact~1. 

\vskip .1in
{\bf Fact 2} \cite{spher}.
\emph{Assume that for every $n\ge m$, the matrix $A_{(n)}\in\mathbb{R}^{m\times n}$ has rank~$m$,
and $\xi_{(n)}$ is a random vector uniformly distributed over the unit ball~$Q$ in~$\mathbb{R}^n$.
Then, as $n\to \infty$, the random vector
$$
\rho_{(n)} = n^{1/2}\bigl( A_{(n)}A^\top_{(n)} \bigr)^{-1/2}A_{(n)}\xi_{(n)}
$$
tends in distribution to the standard Gaussian vector ${\cal N}(0,{\bf I}_m)$, where ${\bf I}_m$ is the identity $m\times m$-matrix.}

Note that for~$n$ fixed, we have
\begin{equation}
\label{rho_rho}
\|\rho_{(n)}\|^2 = n\rho;
\end{equation}
i.e., Facts~2 and~1 characterize the asymptotic distribution of the vector~$\rho_{(n)}$
and exact distribution of its squared norm (normalized by the dimension).

\section{Main Results: Ball-Shaped Sets}
\label{S:ball}

In this section we analyse the two optimization settings formulated in Section~\ref{two_settings} for $Q$ being the $n$-dimensional unit $l_2$-ball.

\subsection{Scalar Optimization}
\label{ssec:scalar}
We consider the scalar case~\eqref{linfun} and discuss first a qualitative result that follows immediately from Fact~1.

Without loss of generality, let $c = (1,\, 0,\,\dots,\, 0)^\top$, so that the function~$g(x)=x_1$ takes its values on the segment $[-1,\, 1]$,
and the true maximum of $g(x)$ on~$Q$ is equal to~$1$ (respectively,~$-1$ for the minimum) and is attained with $x = c$ (respectively, $x=-c$).
Let us compose the random variable
$$
\rho = g^2(\xi),
$$
which is the squared first component~$\xi_1$ of~$\xi$. By Fact~1 with $m=1$ (i.e., $A = c^\top$),
for the probability density function (pdf) of $\rho$ we have
\begin{equation}
\label{scalar_pdf}
f_\rho(x) \,=\, \frac{\Gamma(\frac{n}{2}+1)}{\Gamma(\frac{1}{2})\Gamma(\frac{n+1}{2})}x^{-\frac{1}{2}}(1-x)^{\frac{n-1}{2}} \,\doteq\,
\beta_n\, x^{-\frac{1}{2}}(1-x)^{\frac{n-1}{2}}.
\end{equation}

\begin{figure}[h!]
\centerline{
\includegraphics[width=70mm]{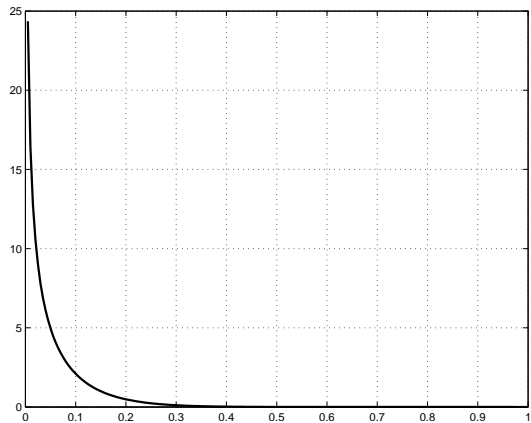}
}
\caption{The probability density function~\eqref{scalar_pdf} for $n=20$.}
\label{fig:pdf_rho}
\end{figure}

Straightforward analysis of this function shows that, as dimension grows, the mass of the distribution tends to concentrate closer the origin,
meaning that the random variable (r.v.)~$\rho$ is likely to take values which are far from the maximum, equal to unity.
To illustrate, Fig.~\ref{fig:pdf_rho} depicts the plot of the pdf~\eqref{scalar_pdf} for $n=20$.

We next formulate the following rigorous result.

\begin{theorem}
\label{th:ball_maxlin}
Let $\xi$ be a random vector uniformly distributed over the unit Euclidean ball~$Q$ and let $g(x)=x_1$, $x\in Q$.
Given $p\in(0,\,1)$ and $\delta\in(0,\,1)$, the minimal sample size~$N_{\min}$
that guarantees, with probability at least~$p$, for the empirical maximum of~$g(x)$ to be at least a~$\delta$-accurate estimate of the true maximum,
is given by
\begin{equation}
\label{Nmin_ball}
N_{\min} = \frac{{\rm ln}(1-p)}{{\rm ln}\Bigl[\frac{1}{2}+\frac{1}{2}I\bigl((1-\delta)^2; \frac{1}{2},\frac{n+1}{2}\bigr)\Bigr]}\,,
\end{equation}
where $I(x; a, b)$ is the regularized incomplete beta function with parameters~$a$ and~$b$.
\end{theorem}

\vskip .1in
Clearly, a correct notation should be $N_{\min} = \lceil \cdot \rceil$, i.e., rounding toward the next integer; we omit it,
but it is implied everywhere in the sequel.

\vskip .1in
\emph{Proof}
We specify sample size~$N$, and let $\left.\xi^{(i)}\right|_1^N$ be a multisample from the uniform distribution on~$Q$;
also introduce the random variable
\begin{equation}
\label{empir_max_ball1}
\eta = \max_{1\leq i \leq N}g(\xi^{(i)}),
\end{equation}
the empirical maximum of the function $g(x) = a^\top x$, $a = (1,\, 0,\,\dots,\, 0)^\top$, from this multisample.
We now estimate the probability ${\sf P}\{\eta>1-\delta\}$.

By Fact~1, the pdf of the r.v.~$\rho=g^2(\xi)$ is given by~\eqref{scalar_pdf},
and its cumulative distribution function (cdf) is known to be referred to as the \emph{regularized incomplete beta function}
$I\bigl(x; \frac{1}{2},\frac{n+1}{2}\bigr)$ with parameters $\frac{1}{2}$ and $\frac{n+1}{2}$, \cite{Wilks}.
Due to the symmetry of the distribution of~$\xi_1$, we have
${\sf P}\{\rho>(1-\delta)^2\} = 2{\sf P}\{\xi_1>1-\delta\}$, so that
${\sf P}\{\xi_1>1-\delta\} = \frac{1}{2}-\frac{1}{2}I\bigl((1-\delta)^2; \frac{1}{2},\frac{n+1}{2}\bigr)$.
Respectively, ${\sf P}\{\xi_1\leq1-\delta\} = \frac{1}{2}+\frac{1}{2}I\bigl((1-\delta)^2; \frac{1}{2},\frac{n+1}{2}\bigr)$ and
${\sf P}\{\eta\leq1-\delta\} = \Bigl[\frac{1}{2}+\frac{1}{2}I\bigl((1-\delta)^2; \frac{1}{2},\frac{n+1}{2}\bigr)\Bigr]^N$, so that finally
\begin{equation}
\label{theProb}
{\sf P}\{\eta>1-\delta\} = 1 - \Bigl[\frac{1}{2}+\frac{1}{2}I\bigl((1-\delta)^2; \frac{1}{2},\frac{n+1}{2}\bigr)\Bigr]^N.
\end{equation}
Letting ${\sf P}\{\eta>1-\delta\} = p$ and inverting the last relation, we arrive at~\eqref{Nmin_ball}. \qed

\vskip .1in
Numerical values of the function $I(x; a, b)$ can be computed via use of the {\sc Matlab} routine {\tt betainc}.
For instance, with modest values $n=10$, $\delta=0.05$, and $p=0.95$, this gives $N_{\min}\approx 8.9\cdot 10^6$, and this quantity grows quickly as
the dimension~$n$ increases.

Since we are interested in small values of~$\delta$, i.e., in $x$ close to unity, a ``closed-form'' lower bound for $N_{\min}$ can be computed as
formulated below.

\begin{corollary}
In the conditions of Theorem~\ref{th:ball_maxlin}
$$
N_{\min} > N_{\rm appr} = \frac{{\rm ln}(1-p)}{{\rm ln}\Bigl[ 1- \tfrac{\beta_n}{n+1}\tfrac{1}{1-\delta}\bigl(2\delta - \delta^2 \bigr)^{(n+1)/2} \Bigr]}\,,
$$
where $\beta_n = \frac{\Gamma(\frac{n}{2}+1)}{\Gamma(\frac{1}{2})\Gamma(\frac{n+1}{2})} = 1/B(\tfrac{1}{2},\tfrac{n+1}{2})$\,.
\end{corollary}

\emph{Proof}
We have
\begin{eqnarray*}
I(x; \tfrac{1}{2}, \tfrac{n+1}{2})
& = &  \beta_n\int_0^x t^{-1/2}(1-t)^{(n-1)/2}{\rm d}t \\
& = &  \beta_n\int_0^1 t^{-1/2}(1-t)^{(n-1)/2}{\rm d}t - \gamma_n\int_x^1 t^{-1/2}(1-t)^{(n-1)/2}{\rm d}t\\
& > &   1 - \beta_n\int_x^1 x^{-1/2}(1-t)^{(n-1)/2}{\rm d}t \qquad \mbox{[ since $t^{-1/2}< x^{-1/2}$ for $x<t < 1$ ]}\\
& = &   1 - \beta_n x^{-1/2}\int_0^{1-x}v^{(n-1)/2}{\rm d}v  \qquad\quad \mbox{[ $v=1-t$ ]}\\
& = &   1 - \beta_n\frac{2}{n+1}x^{-1/2}(1-x)^{(n+1)/2},
\end{eqnarray*}
so that from~\eqref{theProb} we obtain
$$
{\sf P}\{\eta>1-\delta\} \,>\,  1- \Bigl[ 1- \tfrac{\beta_n}{n+1}\tfrac{1}{1-\delta}\bigl(2\delta - \delta^2 \bigr)^{(n+1)/2} \Bigr]^N
$$
and
$$
N_{\rm appr} = \frac{{\rm ln}(1-p)}{{\rm ln}\Bigl[ 1- \tfrac{\beta_n}{n+1}\tfrac{1}{1-\delta}\bigl(2\delta - \delta^2 \bigr)^{(n+1)/2} \Bigr]} < N_{\min}.
$$
Proof is complete. \qed

\vskip .1in
Further simplification of the lower bound can be obtained:
$$
N_{\rm appr} > \widetilde N_{\rm appr} = -\frac{{\rm ln}(1-p)}
{\sqrt{2\pi (n+1)}\tfrac{1}{1-\delta}\bigl(2\delta - \delta^2 \bigr)^{(n+1)/2}}\,.
$$
This is doable by noting that ${\rm ln}(1-\varepsilon)\approx -\varepsilon$ for small $\varepsilon>0$
and using the approximation $B(a,b)\approx\Gamma(a)b^{-a}$ for the beta function with large~$b$.
These lower bounds are quite accurate; for instance, with $n=10$, $\delta=0.05$, and $p=0.95$, we have
$N_{\min}\approx 8.8694\cdot 10^6$, while $N_{\rm appr} \approx 8.7972\cdot 10^6$ and $\widetilde N_{\rm appr} = 8.5998\cdot 10^6$.

\subsection{Multiobjective Optimization}
\label{ssec:multiobj}

Consider a (possibly nonlinear) mapping $g:\mathbb{R}^n\rightarrow \mathbb{R}^m$, $n\gg m>1$; the goal is to characterize the boundary of the image of
a set~$Q\subset\mathbb{R}^n$ under the mapping~$g$. Apart from being of independent interest, this problem emerges in numerous applications.
In particular, if a special part of the boundary,
the \emph{Pareto front}~\cite{Deb} is of interest, we arrive at a multiobjective optimization problem.
Numerous examples (e.g., see~\cite{sinaia}) show that, for~$n$ large, the images of the points sampled randomly uniformly in~$Q$ may happen to fall
deep inside the true image set, giving no reasonable description of the boundary and the Pareto front of $g(Q)$.

We first present a qualitative explanation of this phenomenon by using the setup of Fact~2;
i.e., the set~$Q$ is the unit Euclidean ball and the mappings are linear.

\begin{figure}[h!]
\centerline{
\includegraphics[width=90mm]{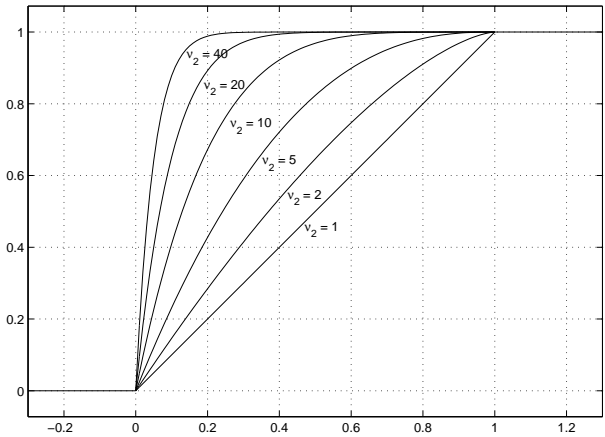}
}
\caption{Cumulative distribution functions ${\cal B}(1, \nu_2)$ for $\nu_2 = 1, 2, 5, 10, 20, 40$.}
\label{fig:beta_distr}
\end{figure}

Since the squared norm of a standard Gaussian vector in $\mathbb{R}^m$ has the $\chi^2$-distribution ${\cal C}(m)$ with~$m$ degrees of freedom~\cite{Wilks},
from Fact~2 and \eqref{rho_rho} we obtain
$$
n\rho \to {\cal C}(m)
$$
in distribution as $n\to\infty$. This is in compliance with the well-known result in the probability theory, namely,
$\nu_2{\cal B}(\nu_1, \nu_2) \rightarrow {\cal C}(\nu_1)$ in distribution as $\nu_2\to\infty$, \cite{Wilks};
here, ${\cal B}(\nu_1, \nu_2)$ stands for the r.v. having the beta distribution with shape parameters $\nu_1,\nu_2$.
For $\nu_1=1$ (i.e., $m=2$, the case most relevant to applications), Fig.~\ref{fig:beta_distr}
depicts the plots of the cumulative distribution functions ${\cal B}(\nu_1, \nu_2)$ (see~\eqref{cdf_rho} below for the explicit formula)
for $\nu_2 = 1, 2, 5, 10, 20, 40$ (i.e., $n = 2,4,10,20,40,80$).

Hence, Fact~2 immediately implies the following important conclusion: Linear transformations essentially change
the nature of the uniform distribution on a ball. Namely, as the dimension of the vector~$\xi$ grows, with the rank of the transformation matrix~$A$ being unaltered, the distribution of the vector~$A\xi$ tends to ``concentrate closer the center'' of the image set.

We now turn to Fact~1 and provide quantitative estimates; to this end, consider the simple case where $m=2$
and the two mappings are linear:
\begin{equation}
\label{lin_mapping}
g_1(x) = c_1^\top x,\quad \|c_1\| = 1, \quad g_2(x) = c_2^\top x,\quad \|c_2\| = 1, \quad c_1^\top c_2=0
\end{equation}
(i.e., $A = \left(
              \begin{array}{c}
                c_1^\top \\
                c_2^\top\\
              \end{array}
            \right)
$ in the notation of Fact~1);
for instance, $c_1, c_2$ may be any two different unit coordinate vectors, so that $g_1(x)=x_i$ and $g_2(x)=x_j$, $i\neq j$, are the two
different components of~$x$. Then the image of~$Q$ is the unit circle centered at the origin.

Introduce now the random variable
\begin{equation}
\label{sq_norm2}
\rho = g_1^2(\xi) + g_2^2(\xi),
\end{equation}
the squared norm of the image of~$\xi\sim {\cal U}(Q)$ under mapping~\eqref{lin_mapping} (i.e., $\rho=\xi_i^2+\xi_j^2$).
Then, by Fact~1 with $m=2$, we have the closed-form expressions for the
cdf~$F_\rho$ and pdf~$f_\rho$ of the r.v.~$\rho$:
\begin{equation}
\label{cdf_rho}
F_\rho(x) =  \left\{
\begin{array}{cl}
           0           & \mbox{~~~for~} x<0, \\
1-(1-x)^{\frac{n}{2}}  & \mbox{~~~for~} 0\leq x\leq1, \\
           1           & \mbox{~~~for~} x > 1;
\end{array}
           \right.
\end{equation}
\begin{equation}
\label{pdf_rho}
f_\rho(x) =  \left\{
\begin{array}{cl}
\frac{n}{2}(1-x)^{\frac{n}{2}-1}
  & \mbox{~~~for~} 0<x<1,  \\
0 & \mbox{~~~otherwise}.
\end{array}
           \right.
\end{equation}

With these in mind, let us evaluate the minimal length~$N$ of the multisample that guarantees a given accuracy with a given probability.
To this end, recall that, given a multisample $\left.\zeta^{(i)}\right|_1^N$ from the scalar cdf~$F_\zeta(x)$ with pdf~$f_\zeta(x)$, the random variable
$$
\eta = \max\{\zeta^{(1)},\dots,\zeta^{(N)}\}
$$
has the cumulative distribution function $F_\eta(x) = F_\zeta^N(x)$ with pdf
$$
f_\eta(x) = F^\prime_\eta(x) = N f_\zeta(x)F_\zeta^{N-1}(x),
$$
which is, in our case \eqref{cdf_rho}--\eqref{pdf_rho} writes
\begin{equation}
\label{eta_pdf}
f_\eta(x) = \frac{Nn}{2}(1-x)^{n/2-1}\Bigl(1-(1-x)^{n/2} \Bigr)^{N-1}.
\end{equation}

We next evaluate several statistics of the r.v.~$\eta=\max\{\rho^{(1)},\dots,\rho^{(N)}\}$.

\bigskip
{\bf Probability:} The theorem below determines the minimal sample size~$N_{\min}$ that guarantees, with high probability,
that a random vector $\xi\sim{\cal U}(Q)$ be mapped close to the boundary of the image set.

\begin{theorem}
\label{th:mapping_prob}
Letting $\xi$ be the random vector uniformly distributed over the unit Euclidean ball $Q\subset\mathbb{R}^n$,
consider the linear mapping $g(\cdot)$ as in~\eqref{lin_mapping}.
Given $p\in(0,\,1)$ and $\delta\in(0,\,1)$, the minimal sample size $N_{\min}$
that guarantees, with probability at least~$p$, that at least one sample be mapped at least~$\delta$-close to the boundary of the image set,
is given by
$$
N_{\min} = \frac{{\rm ln}(1-p)}{{\rm ln}\Bigl(1-(2\delta-\delta^2)^{n/2}\Bigr)}\,.
$$
For small $\delta$ we have
$$
N_{\min} \approx -\frac{{\rm ln}(1-p)}{(2\delta-\delta^2)^{n/2}}\,.
$$
\end{theorem}

\setlength{\baselineskip}{15pt}
\emph{Proof} Let us specify sample size~$N$ and estimate the probability for a sample to be mapped close to the boundary of the image set.
To this end, denote the image of~$\xi\sim{\cal U}(Q)$ under mapping~\eqref{lin_mapping} by $g(\xi) = \bigl(g_1(\xi),\, g_2(\xi)\bigr)^\top$ and
introduce the r.v.
$$
\varkappa = \max\{ \|g^{(1)}\|, \dots, \|g^{(N)}\| \},
$$
the maximum of $\|g(\xi)\|$ over the multisample~$\left.\xi^{(i)}\right|_1^N$.
Also, consider the r.v. $\rho=\rho(\xi)$~\eqref{sq_norm2} for which we have
$$
F_\rho(x^2) = {\sf P}\{\rho\leq x^2\} =  {\sf P}\{\|g(\xi)\|^2\leq x^2\} = 1-(1-x^2)^{\frac{n}{2}}
$$
and the r.v. $\eta = \max\{\rho^{(1)},\dots,\rho^{(N)}\}$,
the maximum of $\rho$ over the multisample~$\left.\xi^{(i)}\right|_1^N$, for which we have
\begin{equation}
\label{eta_cdf}
{\sf P}\{\eta\leq x^2\} = F_\rho^N(x^2) = \Bigl(1-(1-x^2)^{\frac{n}{2}}\Bigr)^N.
\end{equation}

Hence, noting that $\eta=\varkappa^2$, for a small~$\delta>0$ (i.e., letting $x = 1-\delta$),
we see that the probability for at least one sample~$\xi^{(i)}$ to be mapped at least $\delta$-close
to the boundary is equal to
$$
{\sf P}\{\varkappa > 1-\delta\} = 1-\Bigl(1-(2\delta-\delta^2)^{\frac{n}{2}}\Bigr)^N.
$$

Let $p\in(0,\,1)$ be a desired confidence level; then, letting ${\sf P}\{\varkappa > 1-\delta\} = p$ and
inverting the last relation, we obtain the minimal required length of the multisample.
The simple approximation for $N_{\min}$ follows from the fact that ${\rm ln}(1-\varepsilon)\approx -\varepsilon$
for small $\varepsilon>0$.  \hfill $\square$

\begin{figure}[h!]
\centerline{
\includegraphics[width=80mm]{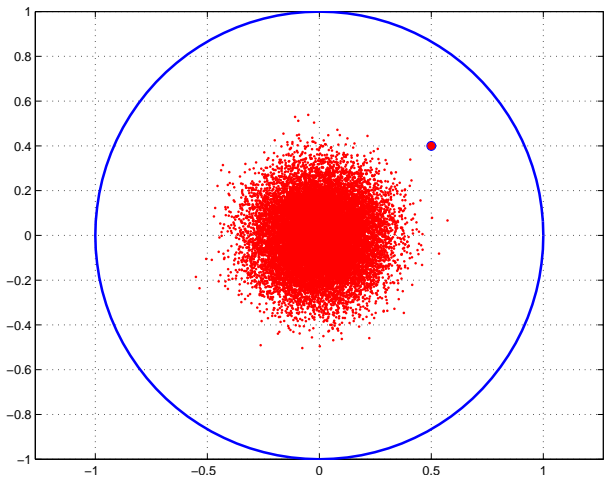}
}
\caption{The 2D image of the $50$-dimensional ball and the result of the Monte Carlo sampling}
\label{fig:twoDmap}
\end{figure}

\vskip .1in
To illustrate, for modest dimension $n=10$, accuracy $\delta=0.05$, and probability $p=0.95$,
one has to generate approximately $N=3.4\!\cdot\!10^5$ random samples to obtain, with probability~$95\%$, a point which is $\delta$-close to the boundary
of the image set.
A sharper illustration of this phenomenon for $n=50$ is given in Fig.~\ref{fig:twoDmap}, which depicts
the images of $N = 100,000$ samples of~$\xi\sim{\cal U}(Q)$ under mapping~\eqref{lin_mapping}.
None of them falls closer than $0.35$ to the boundary of the image set.

\setlength{\baselineskip}{17pt}
Qualitatively, such a behavior can be explained by using geometric considerations and simple
projection-type arguments discussed in~\cite{HK}.

\medskip
{\bf Mode:} The pdf~\eqref{eta_pdf} can be shown to be unimodular, and we find its mode by straightforward differentiating.
Letting $z=1-x$, for the pdf we have
$$
f(z) = \frac{Nn}{2}z^{n/2-1}\Bigl(1-z^{n/2} \Bigr)^{N-1}.
$$
Then $f^\prime(z)=0$ writes
$$
(n/2-1)z^{n/2-2}(1-z^{n/2})^{N-1} = z^{n/2-1}(N-1)(1-z^{n/2})^{N-2}\frac{n}{2}z^{n/2-1}.
$$
Simplifying, we obtain
$$
z^{n/2} = \frac{n-2}{nN-2},
$$
hence,
$$
x_{\max} = 1-\Bigl(\frac{n-2}{nN-2} \Bigr)^{2/n}.
$$

We thus arrive at the following result.

\begin{theorem}
\label{th:mode}
Letting $\xi$ be the random vector uniformly distributed over the unit Euclidean ball $Q\subset\mathbb{R}^n$,
consider the linear mapping $g(\cdot)$ as in~\eqref{lin_mapping} and the random variable
$$
\eta = \max_{i=1,\dots,N}\rho^{(i)},
$$
the empirical maximum of the function $\rho(x)=\|g(x)\|^2$ from the multisample $\left.\xi^{(i)}\right|_1^N$ of size~$N$.
The mode of the distribution of $\eta$  is given by
$$
x_{\max} = 1-\Bigl(\frac{n-2}{nN-2} \Bigr)^{2/n}.
$$
For large~$n$ we have an approximation
$$
x_{\max} \approx 1-\frac{1}{N^{2/n}}.
$$
\end{theorem}

The quantity $x_{\max}$ is seen to be essentially less than unity for large~$n$, even if the sample size~$N$ is huge.
This means the r.v.~$\eta$ takes values far from the boundary of the image.
\begin{figure}[h!]
\centerline{
\includegraphics[width=80mm]{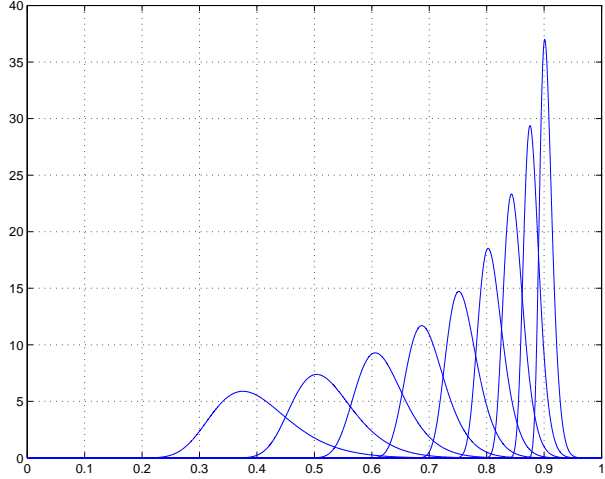}
}
\caption{Probability density functions~\eqref{eta_pdf} for $n=20$ and $N=10^k$, $k = 2,\dots,10$}
\label{fig:eta_pdf}
\end{figure}
For instance, let $n=20$; then, to ``obtain a point in the $0.1$ vicinity of the boundary,''
one has to generate $N\approx10^{10}$ random samples in~$Q$. The family of the pdfs~\eqref{eta_pdf} is plotted in Fig.~\ref{fig:eta_pdf}.

\bigskip
{\bf Expectation:} We now estimate the mathematical expectation~${\sf E}$ of the empirical maximum.

\begin{theorem}
\label{th:expect}
Under the conditions of Theorem~\ref{th:mode} we have
\begin{equation}
\label{expect}
{\sf E}(\eta) = 1 - \frac{2}{n}B\Bigr(\frac{2}{n},N+1\Bigr),
\end{equation}
where $B(\cdot,\cdot)$ is the beta function.
\end{theorem}

\emph{Proof}
If a r.v. $\zeta$ is positive and defined on $D=[0,\;d]$, then the expectation
$$
{\sf E}(\zeta) = \int\limits_0^d\Bigr(1-F_\zeta(x)\Bigl){\rm d}x, 
$$
where $F_\zeta(x)$ is the cdf of~$\zeta$. Hence, having $N$ samples $\xi^{(i)}$ of $\xi\sim{\cal U}(Q)$ and the respective r.v. $\rho=\rho(\xi)$~\eqref{sq_norm2}
with support $[0,\,1]$, for the r.v. $\eta =  \max\{\rho^{(1)},\dots,\rho^{(N)}\}$ we have
$$
{\sf E}(\eta) = \int\limits_0^1\bigl( 1- F_\eta(x)\bigr){\rm d}x,
$$
where $F_\eta(x)$ is given by~\eqref{eta_cdf}. By change of variables $z=(1-x)^{n/2}$, we arrive at \eqref{expect}.\hfill $\square$

For large $n$ and $N$, numerical values of the expectation are close to those observed for the mode;
this is seen from the shape of the pdf~\eqref{eta_pdf} depicted in Fig.~\ref{fig:eta_pdf}.
More formally, having the approximation $B\Bigr(\frac{2}{n},N+1\Bigr)\approx \Gamma(\tfrac{2}{n})(N+1)^{-2/n}$ for large~$N$,
from~\eqref{expect} we obtain
$$
{\sf E}(\eta) \approx 1-\frac{2}{n}\Gamma(\tfrac{2}{n})(N+1)^{-2/n} = 1-\Gamma(\tfrac{2}{n}+1)(N+1)^{-2/n} \approx 1-N^{-2/n}.
$$

For instance, with $n=20$ and $N=10^9$, we have ${\sf E}(\eta) = 0.8802$ for the expectation and $x_{\max} = 0.8754$ for the mode.

\section{Main Results: Box-Shaped Sets}
\label{S:box}

In this section, we consider the scalar optimization problem, however, for box-shaped sets, i.e,, not related to Facts~1 and~2.
We consider the scalar setup described in Section~\ref{two_settings} along with the deterministic approach based on use of regular grids.

\subsection{A Direct Monte Carlo Approach}
\label{ssec:boxMC}

Consider the linear scalar optimization problem for the case where $Q = [-1,\,1]^n$.
Clearly, the results heavily depend on the vector~$c$ in the optimized function $g=c^\top x$;
we consider two extreme cases.

\vskip .1in
{\bf Case 1.} First, let $c=(1,\,0, \dots, 0)^\top$ and consider the empirical maximum
$$
\eta = \max\{ g^{(1)},\dots, g^{(N)} \},
$$
where $g^{(i)}$ is the first component of the random vector $\xi\sim {\cal U}(Q)$.
Specifying~$\delta\in(0,\,1)$, we obtain
$$
{\sf P}\{\eta \leq 1-\delta\} = (1-\delta/2)^N.
$$
This quantity is seen to be independent of the dimension (which is obvious as it is).
Now, specifying a probability level $p\in[0,\, 1]$,  we obtain that the minimal required sample size that guarantees accuracy~$\delta$
with probability $p$ is equal to
$$
N_{\min} = \frac{\ln (1-p)}{\ln (1-\delta/2)}.
$$
For instance, with $p = 0.95$ and $\delta = 0.1$, one has to generate just $59$ points
to obtain a $10\%$-accurate estimate of the maximum  with probability~$95\%$, independently of the dimension.

\vskip .1in
{\bf Case 2.}
Now let $c=(1,\,1, \dots, 1)^\top$; i.e., the optimized function is $g(x) = \sum_i^n x_i$, so that the maximum is attained at
$x=c^\top$ and is equal to $\eta^* = n$. In contrast to Case~1, Monte Carlo sampling exhibits a totally different behavior.
Below, ${\bf Vol}(\cdot)$ stands for the volume of a set.

\begin{theorem}
\label{th:box_diag}
Letting $\xi$ be the random vector uniformly distributed over the unit $l_\infty$-norm ball $Q = [-1,\, 1]^n$,
consider the linear function $g(x) = \sum_i^n x_i$.
Given $p\in(0,\,1)$ and $\delta\in(0,\,1)$, $\delta\leq 1/n$, the minimal sample size~$N_{\min}$
that guarantees, with probability at least~$p$, for the empirical maximum of~$g(x)$ to be at least a~$\delta$-accurate estimate of the true maximum,
is given by
\begin{equation}
\label{Nmin_box_diag}
N_{\min} = \frac{{\rm ln}(1-p)}{{\rm ln}\bigl(1-\frac{n^n \delta^n}{2^n n!}\bigr)}\,.
\end{equation}
For small~$\delta$ and large~$n$ we have
\begin{equation}
\label{stirling}
N_{\min} < -\frac{\sqrt{2\pi n}}{(\delta{\rm e}/2)^n}\,{\rm ln}(1-p).
\end{equation}
\end{theorem}

\emph{Proof}
Let us specify a small $\delta\in (0,\, 1)$ and define
$$
Q_\delta = \{x\in Q\colon \sum_i^n x_i \ge n(1-\delta)\},
$$
so that the maximum of $g(x)$, over $Q\setminus Q_\delta$ is equal to $n(1-\delta)$.
For $\delta\leq 1/n$, the set $Q_\delta$ is seen to be the simplex with $n+1$ vertices at the points $v_0 = (1,\,\dots,\, 1)^\top$ and
$$
v_j = (1,\,\dots, 1,\, \underbrace{1-n\delta}_j, 1, \dots,\,1),  \quad j = 1,\dots,n,
$$
with ${\bf Vol}(Q_\delta) = |\frac{1}{n!}{\rm det}\bigl( v_1-v_0; \dots; v_n-v_0\bigr)| = \delta^n n^n/n!$.
Since ${\bf Vol}(Q) = 2^n$, for $\xi\sim{\cal U}(Q)$  we have
$$
{\sf P}\{\xi\in Q_\delta\} = \frac{\delta^n n^n}{2^n n!} \mbox{~~~~~and~~~~~} {\sf P}\{\xi\in Q\setminus Q_\delta\} = 1-\frac{\delta^n n^n}{2^n n!}\,,
$$
so that
$$
{\sf P}\{\eta > n(1-\delta)\} = 1 - \Bigl(1 - \frac{\delta^n n^n}{2^n n!}\Bigr)^N.
$$
Equating this probability to~$p$ and inverting this relation leads to~\eqref{Nmin_box_diag}.

The lower bound~\eqref{stirling} follows immediately from Stirling's formula and the fact that ${\rm ln}(1-\varepsilon)\approx -\varepsilon$
for small $\varepsilon>0$.\hfill $\square$

\medskip
For  $n=10$ and the same values $\delta=0.1$ and $p=0.95$, we obtain a huge $N_{\min}\approx 1.12\cdot 10^{10}$.
Even for $n=2$, an ``unexpectedly'' large number $N_{\min}\approx 600$ of samples are to be drawn.

\vskip .1in
{\bf $l_1$-norm ball:} The setup of Case~2 is of the same flavor as the one where the set~$Q$ is the unit $l_1$-norm ball, and the optimized function is
$g(x) = c^\top x$ with $c = (1,\, 0,\,\dots, 0)^\top$. We have a result similar to those in Theorems~\ref{th:ball_maxlin} and~\ref{th:box_diag}.

\begin{theorem}
\label{th:l_one_ball}
Letting $\xi$ be the random vector uniformly distributed over the unit $l_1$-norm ball $Q = \{x\in \mathbb{R}^n\colon \sum_{i=1}^n|x_i|\leq 1\}$,
consider the linear function $g(x) = x_1$.
Given $p\in(0,\,1)$ and $\delta\in(0,\,1)$, the minimal sample size~$N_{\min}$
that guarantees, with probability at least~$p$, for the empirical maximum of~$g(x)$ to be at least a~$\delta$-accurate estimate of the true maximum,
is given by
\begin{equation}
\label{Nmin_diamond}
N_{\min} = \frac{{\rm ln}(1-p)}{{\rm ln}\bigl(1-\frac{1}{2}\delta^n\bigr)}\,.
\end{equation}
For small $\delta$ we have
$$
N_{\min} \approx -\frac{{\rm ln}(1-p)}{0.5\,\delta^n}\,.
$$
\end{theorem}

\emph{Proof}
The true maximum of $g(x)$ on~$Q$ is equal to unity; we specify accuracy $\delta\in(0,\,1)$ and consider the set
$$
Q_\delta = \{x\in Q\colon\, x_1\ge 1-\delta\}.
$$
We then have
$$
{\bf Vol}(Q) = \frac{2^n}{n!},\qquad {\bf Vol}(Q_\delta) = \frac{(2\delta)^n}{2\!\cdot \! n!},
$$
so that for $\xi\sim{\cal U}(Q)$ we obtain
$$
{\sf P}\{\xi\in Q\setminus Q_\delta\} = \frac{{\bf Vol}(Q\setminus Q_\delta)}{{\bf Vol}(Q)} = 1-\frac{1}{2}\delta^n,
$$
and the rest of the proof is the same as that of the previous theorem. \hfill $\square$

\vskip .1in

To compare complexity associated with evaluating the optimum of a linear function over the $l_2$-, $l_\infty$-, and $l_1$-balls,
we present a table showing the minimal required number of samples for $\delta=0.05$, $p=0.95$ and various dimensions,
as per formulae~\eqref{Nmin_ball}, \eqref{Nmin_box_diag}, and \eqref{Nmin_diamond}, respectively.
\begin{table}[h!]
\begin{center}
\begin{tabular}{|c|c|c|c|c|c|c|c|}
  \hline
    $n$     &    1  &         2         &          3        &         4         &         5         &         10            &          15        \\
  \hline
   $l_2$    &  119  &        449        &  $1.6\cdot 10^3$  &  $5.7\cdot 10^3$  &  $2\cdot 10^4$    &  $8.9\cdot 10^6$      &  $3.6\cdot 10^9$   \\
 $l_\infty$ &  119  &  $2.4\cdot 10^3$  &  $4.3\cdot 10^4$  &  $7.2\cdot 10^5$  &  $1.2\cdot 10^7$  &  $1.1\cdot 10^{13}$   &  $10^{19}$ \\
   $l_1$    &  119  &  $2.4\cdot 10^3$  &  $4.8\cdot 10^4$  &  $9.6\cdot 10^5$  &  $1.9\cdot 10^7$  &  $6.1\cdot 10^{13}$   &  $2\cdot 10^{20}$  \\
  \hline
\end{tabular}
\label{T:l_p_boxes}
\vspace{.1in}
\caption{$l_p$-balls: Minimal required number of samples for $\delta=0.05$ and $p=0.95$.}
\end{center}
\end{table}

These results are in consistence with intuition, since the $l_1$-norm ball is ``closer'' in shape to the ``worst-case'' conic set, 
while the $l_\infty$-norm ball with $c=(1,\,1, \dots, 1)^\top$ ``takes an intermediate position'' between $l_2$ and $l_1$ (obviously, closer to~$l_1$).

\subsection{Deterministic Grids}
\label{ssec:lptau}

In this section we briefly discuss a natural alternative to random sampling the set~$Q$. A belief is that use of various deterministic grids
might outperform straightforward Monte Carlo.
Again, we consider~$Q$ being the unit box $[-1,\, 1]^n$, and the scalar function to be optimized is $g(x) = \sum_{i=1}^n x_i$,
so that the maximum is equal to~$n$. We show that, even in such a simple setting, deterministic grids happen to be computationally intensive
in high dimensions.

\vskip .1in
{\bf Uniform grid:} Consider a positive integer $M>1$ and the uniform mesh~$\cal M$ on~$Q$, with cell-size~$\Delta=2/(M+1)$;
the mesh is assumed not to cover the boundary of~$Q$. The total amount of points in the mesh is ${\tt card} \,{\cal M}=M^n$ and the maximum
of $g(x)$ over~$\cal M$ is equal to $g_{\cal M} = n(1-\Delta)$. To guarantee relative accuracy~$\delta$ of approximation,
i.e., $g_{\cal M} = n\!\cdot\!(1-\delta)$, one needs cell-size to be $\Delta = \delta$, hence, the overall number of mesh points is equal to
\begin{equation}
\label{unif_grid}
{\tt card} \,{\cal M} = \Bigl( \frac{2}{\delta} - 1 \Bigr)^n\,,
\end{equation}
which amounts to a huge $M \approx 6.13\cdot 10^{12}$ for modest $n=10$ and $\delta=0.1$. Interestingly, to obtain the same accuracy
\emph{with probability $p=0.99$}, Monte Carlo requires ``just'' $N_{\min}=1.7\cdot 10^{10}$ samples!

\vskip .1in
{\bf Sobol sequences:} Another type of grids that can be arranged over boxes are \emph{low-discrepancy} point sets
or \emph{quasi-random} sequences~\cite{Niederreiter}. In practice, they share many properties of
\emph{pseudorandom} numbers, e.g., such as those produced by the {\tt rand} routine in {\sc Matlab}.

Among the variety of quasi-random sequences, so-called \emph{$LP_\tau$ sequences} introduced by I.M. Sobol in 1967,~\cite{Sobol-1967-4-UCMaMP}
(also see~\cite{PSI,StatnikovEtAl-JOTA} for recent developments) are widely used in various application areas.
This sophisticated mechanism heavily exploits the box shape of the set; it is much more efficient than purely deterministic uniform grids
and may outperform straightforward Monte Carlo.

In the experiments, we considered the function $g(x) = \sum_{i=1}^n x_i$ defined on $Q=[-1,\, 1]^n$ and computed its maximum value
over the points of an $LP_\tau$ sequence of length $N=10^6$ for various dimensions; this was performed by using the {\sc Matlab} routine
{\tt sobolset}. The corresponding results are given in the row ``$LP_\tau$'' of Table~2.
The row ``Monte Carlo'' presents empirical maxima obtained by using Monte Carlo sampling with the same sample size~$N$
(averaged over $100$ realizations), and the row ``Uniform grid'' relates to the uniform mesh of cardinality~$N$.
\begin{table}[h!]
\begin{center}
\begin{tabular}{|c|c|c|c|c|c|c|c|}
  \hline
  $n$; true max    &     2    &     3    &     4    &     5    &    10    &    15    &    20     \\
  \hline
   Uniform grid    &  1.9960  &  2.9406  &  3.7576  &  4.4118  &  6.0000  &  7.5000  &  6.6667   \\
 $LP_\tau$         &  1.9999  &  2.9792  &  3.8373  &  4.6844  &  7.9330  &  10.2542 &  10.9470  \\
   Monte Carlo     &  1.9974  &  2.9676  &  3.8731  &  4.6981  &  7.8473  &  10.0796 &  11.8560  \\
  \hline
\end{tabular}
\label{T:box}
\vspace{.1in}
\caption{$l_\infty$-box: Empirical maxima for the three methods with $N=10^6$.}
\end{center}
\end{table}

It is seen that the uniform mesh exhibits a very poor relative performance as dimension grows, while the results of $LP_\tau$ and Monte Carlo approaches
are much better and similar to each other. Clearly, the absolute values are very far from the true maxima equal to~$n$, since $N=10^6$
samples are not sufficient to obtain reasonable accuracy for dimensions $n> 5$.

Instead of computing the sample size for fixed values of the accuracy~$\delta$ as in Table~1,
here we fix the sample size and compute the empirical maxima.
The reason for such an ``inversion'' is that, given~$\delta$ and the specific linear function $g(x)$, it is not quite clear how to estimate the required length of the $LP_\tau$-sequence.
To overcome this problem, one might fix a reasonable value of accuracy, compute the minimal sample size~$N_{\min}$ required for Monte Carlo,
run the $LP_\tau$-approach with this length, and compare the results. However, for large dimensions~$n$, the values of~$N_{\min}$ are huge,
leading to very large computation times or even memory overflow.

The values obtained for the uniform mesh were computed by inverting relation~\eqref{unif_grid} and using the actual grid with cardinality
${\tt card} \,{\cal M} = \lceil N^{1/n}\rceil^n \ge N$, so the quantities presented in row ``Uniform grid'' of Table~2 are overly optimistic.

More importantly, the routine {\tt sobolset} has several parameters, such as $s={\tt skip}$ (choice of the initial point in the sequence)
and $\ell={\tt leap}$ (selection of every $\ell$th point).
These play the similar role as the seed value in {\sc Matlab} pseudorandom number generators does,
and may happen to be crucial for the quality of the resulting $LP_\tau$ sequence. In the experiments, we used the default values $s=\ell=1$
(i.e., the first~$N$ points), though for different values of $s, \ell$, the respective estimates of the empirical maximum may differ a lot.
For instance, adopting $\ell=133$, for $n=5$ we obtain a much better estimate $4.8653$,
while taking (intentionally bad) $\ell=128$ leads to very poor $2.8734$.

Finally, note that applications of uniform grids and $LP_\tau$ sequences are limited to box-shaped sets. Sets different from boxes can in principle
be embedded in tight enclosing boxes with subsequent use of rejection techniques; however, the rejection rate usually grows dramatically
as the dimension increases.

\section{Conclusions}
The main contribution of the paper is a rigorous explanation of the reason why does
a direct Monte Carlo approach show itself inefficient in high-dimensional optimization problems when
estimating the maximum value of a function from a random sample in the domain of definition. First,
attention was paid to linear functions and ball-shaped sets; using known results on the uniform distribution
over the ball, we characterized the accuracy of the estimates obtained via a specific random variable
associated with the function value. Also, a multiobjective optimization setup was discussed. The results
obtained testify to a dramatic growth of computational complexity (required number of samples) as the
dimension of the ball increases. Same flavor results are obtained for box-shaped sets; these also include
analysis of deterministic grids.


The authors are indebted to an anonymous reviewer for the critical comments that led to a tangible progress in the presentation of the results.

\end{document}